\documentclass[12pt,a4paper]{article}
\input{epsf}

\usepackage[T1]{fontenc}
\usepackage{amsmath,amssymb,amsfonts,graphicx,amsthm,mathrsfs}
\usepackage{cite}
\usepackage[T1]{fontenc}
\usepackage[utf8]{inputenc}
\usepackage{tikz}
\usepackage{float}
\usepackage{caption}
\usepackage{authblk}
\usepackage[colorlinks]{hyperref}

\newtheorem{theorem}{Theorem}

\newtheorem{introtheorem}{Theorem}

\newtheorem*{thm*}{Theorem \ignorespaces}

\newtheorem{definition}{Definition}

\newtheorem*{conda}{Condition $\mathbf{A}$}
\newtheorem*{condb}{Condition $\mathbf{B}$}

\theoremstyle{definition}

\newcommand\restr[2]{{% we make the whole thing an ordinary symbol
\left.\kern-\nulldelimiterspace % automatically resize the bar with \right
#1 % the function
\littletaller % pretend it's a little taller at normal size
\right|_{#2} % this is the delimiter
}}

\newcommand{\littletaller}{\mathchoice{\vphantom{\big|}}{}{}{}}

\usepackage{enumerate}
\allowdisplaybreaks
\usepackage{hyperref}
\usepackage{stackrel}
\DeclareMathOperator{\argmin}{argmin}
\DeclareMathOperator{\tr}{tr}
\DeclareMathOperator{\diag}{diag}
\usepackage{array}

\date{}
\begin{document}

\title{Mean Field Variational Bayesian Inference and Statistical Mechanics of Gaussian Mixture Model }
\author[1]{Alireza Bahraini\thanks{bahraini@sharif.edu}}
\author[1]{Saeed Sadeghi\thanks{saeedsadeghi91@gmail.com}}
\affil[1]{ Department of Mathematical Sciences, Sharif University of Technology, P.O.Box 11155-9415, Tehran, Iran.}
\renewcommand*{\Authands}{, }
\maketitle
\begin{abstract} 
One of the main modeling in many data science applications is the Gaussian Mixture Model (GMM), and Mean Field Variational Bayesian Inference (MFVBI) is classically used for approximate fast computation. 
In this paper, our aim is to lay a mathematical foundation for a rigorous analysis of the MFVBI applied to the GMM. 
Several fundamental key concepts from statistical mechanics surge naturally throughout our process. 
It turns out that GMM can be considered as a generalization of Curie-Weiss model in statistical mechanics. The standard quantities like 
partition function, Legendre transform and free energy come into operation. The initial system of equation \eqref{introp} reduces to a simpler \eqref{bonyi}. We introduce a temperature parameter in order to accommodate a phase transition 
phenomena which can finally guarantee the accuracy of the solutions to MFVBI.

\end{abstract}
\textbf{keywords:} gaussian mixture model, geodesic convexity, optimal transport, mean field variational inference (MFVI), partition function

\section{Introduction}
Mean Field Variational Bayesian Inference (MFVBI) is a method for approximating the posterior probability density in Bayesian Statistics. 
In despite of its popularity and successful performance in practice, it suffers from lack of a rigorous uncertainty quantification. 

In order to overcome challenges of Bayesian inference for computing posterior distributions two principal approaches have been proposed both having their roots in statistical mechanics. The Markov chain Mont Carlo (MCMC) method \cite{metro,hasti} which is based on generating consistent samples from the posterior distribution. This method which has been widely applied as a standard tool for many Bayesian inference problems since its development, is also theoretically well-founded and rather simple to analyze. However the MCMC computational cost is very high when applied to large data models \cite{blei}.

Variational Bayesian Inference (VBI) approach has been successfully employed as an alternative method for approximate computation of posterior distributions \cite{waterhouse,attias}. This method which can be applied for complex models of large data size approximates the given posterior distribution throughout a simpler family, named variational distributions. For Mean Field Variational Bayesian Inference (MFVBI) the variational space comprises the space of factorized probability distributions. This technique has been applied to various problems ranging from graphical models to large-scale document analysis, computational neuroscience and computer vision \cite{blei}. In despite of all successful application VBI approach suffers from lack of theoretical support.

In fact the well-known mean field equations discussed in statistical mechanics are a particular case and source of inspiration for MFVBI in its most general form in statistical learning. Our study in this paper reveals a deep connection between GMM and the co-called Curie-Weiss model which is a basic paradigm with in spin glass theories (See \cite{tal} for a rigorous analysis of Curie-Weiss model).
This connection helps transforming the mean field equations associated with GMM into a system of equations of significantly simpler nature. The new perspective over the MFVBI leads to several immediate theoretical consequences such as providing a precise benchmark
for verifying whether the solutions are reliable or not as well as suggesting a modification of the system by introducing a temperature parameter whose adjustment leads to asymptotically correct approximate parameters. 

To be more specific assume that we have a data set $\{x_i\}_{i=1} ^N$ which is generated by a mixture of Gaussian distributions of the form 
\begin{equation}\label{ptrue1}
\tilde{p}_{N}(x)=\sum_{k=1} ^K \tilde{\pi}_k \ \mathcal{N}(x|\tilde{\mu}_k, \tilde{\Lambda} ^{-1}_k),
\end{equation}
where $ \{\tilde{\pi}_k\}_{k=1} ^K$, $\{\tilde{\mu}_k\}_{k=1} ^K$ and $\{\tilde{\Lambda}_k\} _{k=1} ^K$ denote, respectively, the weights, means and the covariance matrices of the $K$ sub-populations in the model. 

The factorization structure in terms of which MFVBI is fabricated for GMM is of the continuous-discrete type $\mathcal{H}=M\times \mathscr{Z}$ where $M$ is a Riemannian manifold and $\mathscr{Z}$ is a finite discrete set. In the GMM the manifold $M$ is the underlying space for continuous parameters (weights, means and covariance matrices) and $\mathscr{Z}$ comprises the different ways we can partition the data into $K$ classes. The probability distribution can be represented in the form 
\begin{equation}\label{eqmu}
\frac{d\mu}{d\omega_g}:=\frac{e^{-\lambda \Phi}}{Z},
\end{equation}
where $\lambda$ is a parameter that grows to infinity with $N$ the number of the data set $\{x_i\}_{i=1} ^N$. Also $\Phi: M\times \mathscr{Z}\rightarrow \mathbb{R}$ is $C$-convex map over each of the sheets $M\times \{i\}$ for $i\in \mathscr{Z}$ where $C$ is a constant independent of $N$ and $i$. In fact we represent the standard GMM as $\lambda_0 \Phi$ and the parameter $\lambda $
equals $\beta \lambda_0$. The parameter $\beta$ is included to play a role similar to temperature parameter in statistical mechanics which is responsible for phase transition phenomenon. It is worthy to recall that in all models from statistical mechanic including Curie -Weisse model, the naive mean field approximation is consistent only for appropriate range of temperature parameter which is characterized by phase transition phenomenon. A similar result holds true for GMM. We will observe that by calibrating $\beta$ one can ensure asymptotically rigorous approximation of the parameters.
So $\beta $ is bounded while $\lambda_0$
tends to infinity with $N$.

The space of factorized probability measures denoted by $\mathscr{A}$ is given by
\[
\mathscr{A}:=\mathcal{P}_2 (M)\times \mathcal{P}(\mathscr{Z}),
\]
where $\mathcal{P}_2(M)$ represents the length space $(\mathcal{P}(M),W_2)$ consisting of the space of Borel probability measures on $M$ equipped with the so-called Wasserstein metric of order 2, denoted by $W_2$, and $\mathcal{P} (\mathscr{Z})$ is the simplex of probability measures in $\mathbb{R}^{|\mathscr{Z}|}$. As is well-known, the variational
Bayesian method consists of approximating a given probability measure $\mu\in\mathcal{P}( M\times \mathscr{Z})$ by an element in $\mathscr{A}$ optimizing the Kullback-Leibler distance 
\begin{equation}\label{mfvbie}
\nu_0:=\arg \min_{\nu\in \mathscr{A}}D_{KL} (\nu\|\mu).
\end{equation}

The basic question which has been remained unanswered is that are there any relation between mode and moments of $\mu$ and $\nu_0$? In the present paper we are concerned about the question regarding the modes and we investigate under what circumstances 
the mode of $\mu$ and $\nu_0$ coincide asymptotically. 
In our previous paper \cite{bah} we used the Lott--Villani--Sturm \cite{vilon} theory of optimal transport to investigate the convexity of the Kullback-Leibler functional $D_{KL} (.\|\mu):\mathscr{A}\rightarrow \mathbb{R}$ corresponding to the GMM. 
Here we start by observing in section \ref{stion3} that, if we replace $\mathscr{Z}$ by a subset $\mathscr{Z}_0$ with the two following properties then all the maps $(-\Phi)|_{M\times \{i\}},$ for all $i\in \mathscr{Z}_0$ will enjoy geodesic convexity: 

\begin{itemize}
\item[i)] The minimum ratio of the elements in each of the partition components has a fixed positive lower bound independent of $N$, 

and

\item[ii)] the means and the precision matrices lie in a convex bounded subset with respect to their geometry.
\end{itemize}

In particular we can conclude that for any $i\in \mathscr{Z}_0$, the restriction $(-\Phi)|_{M\times \{i\}} $ admits a unique minimum
at a point $(m_i,i)\in M\times\{i\}$. 
Now assume that $\nu_0:=\mu^1\times \mu^2$ denotes the solution to the mean field variational equation \eqref{mfvbie}. Then by variational method $\mu^1$ and $\mu^2$ satisfy
\begin{equation}\label{introp}
\begin{cases}
\log \mu^2 (i) &=-\lambda\int_M \Phi (\xi,i) \frac{d\mu^1 }{d\omega_g} d\omega_g-\log Z_{2},\\
\log \frac{d\mu^1 }{d\omega_g}(\xi) &=-\lambda\sum_i \Phi (\xi,i) \mu^{2} (i)-\log Z_{1},
\end{cases}
\end{equation}
where $Z_1$ and $Z_2$ are the normalization constants. From the second equation, it can be seen that $\frac{d\mu^1 _{} }{d\omega_g}(\xi)$ takes the following form 
\[
\frac{d\mu^1 _{} }{d\omega_g}(\xi) =\frac{e^{-\lambda\sum_i u_i \Phi (\xi,i)}}{Z_1},
\]
where $u_i$'s are positive numbers satisfying $\sum_i u_i =1$, and therefore $ \frac{d\mu^1 _{} }{d\omega_g}(\xi) $ is also $-\log$-convex and admits a unique maximum. In Theorem \ref{theoo2} of section \ref{sec4}, it will be shown that the unique maximum of $ \frac{d\mu^1 _{} }{d\omega_g}(\xi) $ is a critical point of the following marginal partition function
\begin{equation}
Z(\xi):= \sum_i e^{-\lambda \Phi (\xi,i)+ R_i },
\end{equation}
where $R_i = O(\sqrt{\lambda})$. 

In section \ref{zedxi} we will provide some insight into how the summation $Z(\xi)$ converges in the weak sense towards an effective partition function denoted by $Z_{eff}$ defined on the space of $K\times K$ Markov matrices. This will motivate us to call the map $-\log Z_{eff}$, the free energy of the GMM. 
More precisely in order to investigate the asymptotic behaviour of the above summation we divide the collection $\mathscr{Z}_0$ into some blocks $\mathscr{Z}_0=\cup_{\hat{A}_d\in \hat{\mathcal{A}}_d}\mathscr{Z}_{\hat{A}_d}$
in such a way that for a fixed $\hat{A}_d \in \hat{\mathcal{A}}_d$, the maps $\Phi|_{{M}\times \{i\}}$ for $i\in \mathscr{Z}_{\hat{A}_d}$ have a negligible oscillation around their average in prbability with respect to the selection of the data $\{x_i\}_{i=1} ^N$.
This can be fulfilled by employing U-statistics as a generalization of central limit theorem. Furthermore one can establish a natural correspondence
\[
\hat{}:\mathcal{A}_d\rightarrow \hat{\mathcal{A}}_d, \hspace{1cm} \hat{} ( A)= \hat{A}
\]
between $\hat{\mathcal{A}}_d$ and a lattice $\mathcal{A}_d$ in the whole space of Markov $K\times K$ matrices represented by $\mathcal{A}$. 
We denote the average $\langle \Phi(i, \xi)\rangle$ for $i\in \mathscr{Z}_{\hat{A}_d}$ by $\hat{\phi}_{A_d} (\xi)$.
\[
\hat{\phi}_{A_d} (\xi)=\langle \Phi(i, \xi)\rangle
\]

This map depends linearly on $A_d$ and can be naturally extended to the whole space $\mathcal{A}$ as a linear map.
The latter is indicated by $\hat{\phi}_A (\xi)$ and is defined by \eqref{fiku}. 
Now let 
\[
\mathfrak{M}:\mathcal{A}\rightarrow M,
\]
denote the map defined by 
\[
\mathfrak{M}(A):=\argmin _{\xi\in M} \hat{\phi}_A (\xi),
\]
and let $\hat{\Phi}: \mathcal{A}\rightarrow \mathbb{R},$ be defined as the Legendre transform of $\hat{\phi}_A (\xi)$:
\begin{equation}\label{phihatag1}
\hat{\Phi}(A):=\min_{\xi\in M}\hat{\phi}_A (\xi)=\hat{\phi}_A (\mathfrak{M} (A)).
\end{equation}
Then $Z_{eff}$ can be described by the following relation
\begin{equation}
Z_{eff} (A) =e^{-\lambda \big (\hat{\Phi} (A)+\psi (A)\big )+ O(\sqrt{ \lambda})},
\end{equation}
Moreover $\psi:\mathcal{A}\rightarrow \mathbb{R}$ is included based on counting the iterations that arise for each subclass $\mathscr{Z}_{\hat{A}_d}$. (See \eqref{sydef})

The computation of partition function is of high importance in statistical mechanics and we postpone addressing it for a future paper. Although our computation of partition function is non rigorous in the form stated here meanwhile the following theorems underlines its role in the behavior of the solutions to the MFVI equations. 

In theorem A we show that the unique maximum of $\frac{d\mu^1}{d\omega_g} (\xi)$ corresponds to a critical point of the free energy $\log Z_{eff}$:

\begin{introtheorem}
\label{thm:A}
If the maximum of the right hand side of \eqref{minmu}
\begin{equation}\label{minmu22}
\log \left ( \frac{d\mu ^1 _{}}{d\omega_g} (\xi)\right )= -\lambda \sum_i \Phi (\xi,i)e^{-\lambda \Phi (m,i)+ R_i-\log Z_2}-\log Z_1,
\end{equation}
with respect to $\xi$ occurs at a point which corresponds to an interior point of $\mathcal {A}$ under the map $\mathfrak{M}$, then in large $\lambda$ limit this point converges towards a critical points of $A\rightarrow (\hat{\phi}_A (\mathfrak{M} (A) +\psi (A)) $. Similarly if the maximum occurs at a corresponding interior point of any of sub-complexes of $\mathcal{A}$ it will coincide with a critical point of the restriction of $A\rightarrow (\hat{\phi}_A (\mathfrak{M} (A) +\psi (A)) $ to that sub-complex. 
\end{introtheorem}

If we define 
\[
\mathfrak{J}:\mathscr{Z}\rightarrow \mathcal{A},
\]
by
\[
\mathfrak{J}(i)= A_d, \hspace{1cm} \text{for } i\in \mathscr{Z}_{\hat{A}_d}.
\]
(See \eqref{j1} for the definition of $\hat{}$ ) then for large values of $\lambda$ the solutions to the system of equations MFVBI \eqref{introp} will approach the following Dirichlet distributions: 

\begin{equation}\label{diracint}
\frac{d\mu^1}{d\omega_g}\simeq \delta (\xi-m), \hspace{0.5cm} \mu^2\simeq \mathfrak{J}^* \left [\frac{|\mathscr{Z}_{\hat{A}_d}|}{N} \delta (A-A_{d_1}) \right ]
\end{equation}

If we define $\mathfrak{N}: M\rightarrow \mathcal{A}$ as the Legendre transform of $\xi\rightarrow \hat{\phi}_A (\xi ) +\hat{\psi} (A) $
\begin{equation}
\mathfrak{N} (\xi ):=\argmin_A \big ( \hat{\phi}_A (\xi ) +\hat{\psi} (A) \big ) 
\end{equation}

The following theorem characterizes the points where the above Dirac measures concentrate around.

\begin{introtheorem}
\label{thm:B}
The two Dirac measures given by the equation \eqref{diracint} constitute a solution to the MFVBA equations \eqref{introp}
iff $B$ and $m$ satify the following system of equation
\begin{equation} \label{bonyi}
\begin{cases} \mathfrak{M} (B)=m\\
\mathfrak{N}(m)=B
\end{cases}
\end{equation}
\end{introtheorem}

\begin{introtheorem}
\label{thm:C}
The solution to the system of equations \eqref{bonyi} gives rise to the true parameters $\{\tilde{\mu}_k, \tilde{\Lambda_k}, \tilde{\pi}_k\}_{k=1,\dots,K}$
if and only if $B$ coincides with one of the vertices of the simplex $(\Delta)^K$ where $\Delta=\{(x_1,\dots,x_{K})|\sum x_i =1\}$. 
\end{introtheorem}

We have introduced a temperature parameter $\beta$ which permits us to concld the following theorem

\begin{introtheorem}
\label{thm:D}
By adjusting the temperature parameter $\beta$ one can ensure that the solution to MFVI approaches the correct value of the parameters for $N$ large enough.
\end{introtheorem}
\section{Notation for GMM}\label{secgmm}
The probability factorization structure in terms of which MFVI is fabricated for this model turns out to be of the continuous-discrete type $\mathcal{H}=M\times \mathscr{Z}$ where $M$ is a Riemannian manifold, and $\mathscr{Z}$ is a finite discrete set. More precisely the manifold $M$ is defined by 
\begin{equation}\label{defm}
M:= \tilde{\mathcal{PV}}^K\times C_K ,
\end{equation}
where $C_K := \{(\pi_1,\dots,\pi _K)\ |\ \pi_k \geq 0,\ \sum_{k=1}^K \pi_k =1\}$ and 
\begin{equation}
\tilde{\mathcal{P}\mathcal{V}}=\mathbb{R}^p \times S_{++} ^P=\left\{(\mu ,\Lambda) \ |\ \ \mu\in\mathbb{R}^P ,\ \Lambda \in S^P _{++} \right\},
\end{equation}
where $S^P _{++}$ denotes the space of symmetric positive definite matrices of dimension $P$.

Therefore the model consists of $K-$components mixture of $P-$dimensional multivariate normals with unknown components involved in $\mathscr{Z}$, means $\mu_1,\dots,\mu_K$, precision matrices $\Lambda_1,\dots,\Lambda_K$, and weights $\pi_1,\dots,\pi_K$ representing the probabilities of the components $1,\dots,K$ respectively. The parameter $N$ is the number of data points, $x$ is a $P-$dimensional vector and $x_n$ is the $n$th observed $P-$dimensional data point. The finite set $\mathscr{Z}$ denotes the finite space
\[
\mathscr{Z}= \{1,\dots,K\}^N ,
\]
and for any element $\zeta=(\zeta_1,\dots,\zeta_N)\in \mathscr{Z},$ there is an associated a sequence 
\begin{align}\label{defz}
z:=\left\{z_{ik}\ |\ \ 1\leq i\leq N, \ 1\leq k\leq K\right\},
\end{align}
where 
\[
z_{ik}=\begin{cases}1&\text{if }\zeta_ i=k,\\ 
0&\text{otherwise}.\end{cases}
\]

So the data generating process is detailed as follows:
\begin{align*}
P_N(x|\mu , \pi , \Lambda)&= \prod _{n=1} ^N P_N (x_n | z_n , \mu , \Lambda)\prod_{k=1 } ^K P_N(z_{nk}| \pi _k),
\\
\log P_N (x_n | z_n , \mu , \Lambda )& = \sum _{n=1} ^N z_{nk} \log \phi_k (x_n) +\tilde{C},
\\
\log \phi_k (x) &= -\frac{1}{2} (x- \mu_k )^T \Lambda _k (x - \mu _k ) + \frac{1}{2}\log |\Lambda _k | +\tilde{C},
\\
\log P_N(z_{nk}|\pi _k )&= \sum _{k=1} ^K z_{nk}\log \pi_k +\tilde{C},
\\
P_N(z_{nk} =1)&=\pi _{k},\quad p(x|\pi, \mu , \Lambda , z) =\prod_{n=1}^N \prod_{k=1}^K \mathcal{N}(x_n |\mu_k , \Lambda _k ^{-1})^{z_{nk}} ,
\end{align*}
and
\begin{align}\label{log}
\log P_N (z, \mu, \pi, \Lambda | x ) &= \sum\limits_{n=1} ^N \sum\limits_{k=1} ^K \left(z_{nk} (\log \pi _k )-
z_{nk}\frac{1}{2} (x_n - \mu _k)^T \Lambda _k (x_n - \mu _k) + \frac{1}{2}z_{nk} \log |\Lambda _k |\right)
\notag\\
&\ +\sum\limits _{k=1} ^K \log p (\mu _k )+ \sum\limits _{k=1} ^K \log p (\Lambda _k ) +\log p (\pi) +\tilde{C}.
\end{align}
Here the respected prior models are considered as follows: a multivariate normal prior for $\mu_k$, a Wishart prior for $\Lambda _k$, and a Dirichlet prior for $\pi$. (See \cite{mj} for more details.)
\\
The standard variational assumption on this mixture model is that 
\begin{equation}
q(\mu , \pi , \Lambda , z ) =\prod _{k=1} ^K q(\mu _k) q(\Lambda _k) q(\pi _k) \prod _{n=1} ^N q(z_n).
\end{equation}

We take $\Phi := -\log P_N$. 
%%%%%%%%%%%%%%%%%%%%%%

\section{Convexity of Gaussian Mixture Model}\label{stion3}
According to \cite{bah} and the references therein the well-known Mean Field Variational Bayesian Approximation Inference (MFVBI), 
corresponds to the case where the Polish space $\mathcal{H}$ is factorized as $\mathcal{H}=\prod_{i=1} ^K \mathcal{H}_i$ into a product of Polish subspaces $\mathcal{H}_i \subset \mathcal{H} $ for $i=1,\dots,K$. 
Let $\mathcal{P}(\mathcal{H}_i)$ denote the space of Borel probability measure on $\mathcal{H}_i$. We set 
\begin{equation}\label{fac}
\mathcal{A}:=\prod_{i=1} ^K \mathcal{P}(\mathcal{H}_i).
\end{equation}

The MFVBI consists of the following minimization problem 
\begin{equation}\label{2}
\arg\min _{\nu \in \mathcal{A}} D_{KL} (\nu\| \mu ).
\end{equation}

Applying variational method, one can describe the solution to the above problem by the following system of equations
\begin{equation}\label{sol}
\log \left(\frac{d\nu_i }{d\omega_i} (x_i) \right) = \mathbb{E}^{\nu _{\backslash i} } \left(\log \frac{d\mu}{d\omega}\right).
\end{equation}
where $\omega =\prod _i \omega_i \in\mathcal {P}(\mathcal{H})$ is a given fixed measure for $\omega_i \in \mathcal{P}(\mathcal{H}_i)$ and $\nu_{\backslash i}=\prod_{j\neq i} \nu _j$.

In reference \cite{bah} we have studied convexity of the functional $\nu\rightarrow D_{KL} (\nu||\mu)$ in two cases: the case where $\mathcal{H}_i $'s are Riemannian manifolds and the case where the above factorization has a hybrid discrete-continuous form. We then introduced a correction to GMM to turn it into a $-\log$ convex distribution. In this section we aim to show that instead of making a correction, the convexity of GMM is in fact satisfied within two natural circumstances expressed in conditions $\mathbf{A}$ and $\mathbf{B}$ below.
Briefly, condition $\mathbf{A}$ concerns the boundedness of parameters and condition $\mathbf{B}$ ensures that in each class there exists enough samples. 
\subsection{Constrains for Convexity }

We would like to demonstrate in this section is that under the two conditions $\mathbf{A}$ and $\mathbf{B}$ below, 
\[
-\log P_N= \lambda \Phi, 
\]
where $\Phi$ is a $C-$convex function for some $C>0,$ independent of $N$ and $\lambda$ grows to infinity as $N$ approaches infinity.

The first condition is about the norm of the precision matrices $\Lambda_k,$ and the means $\mu_k$ as described in the following.
\begin{conda}
We assume that the means and the precision matrices have bounded coefficients and lie in a convex subset with respect to their geometry. This condition in practice can be fulfilled by modifying priors distributions.
\end{conda}
More precisely, we consider the space
\begin{equation}
\tilde{\mathcal{P}\mathcal{V}}=\mathbb{R}^p \times S_{++} ^P=\left\{(\mu ,\Lambda) \ |\ \ \mu\in\mathbb{R}^P ,\ \Lambda \in S^P _{++} \right\},
\end{equation}
and for some positive real number $R$, we define the subset $\mathcal{PV}$ of $\tilde{\mathcal{PV}}$ by 
\begin{equation}\label{mv}
\mathcal{PV}=\left\{(\mu ,\Lambda) \in \tilde{\mathcal{P}\mathcal{V}}\ \left|\ \ |\mu|<R, \ d_{RF}(\Lambda, Id)< R \right.\right\}.
\end{equation}
Here $Id $ is the identity matrix and $d_{RF}(.,.)$ denots the distance with respect to Rao-Fisher metric over $S^P_{++}$ (See \cite{bah}).
Let 
\begin{equation}\label{defm}
M:= \mathcal{PV}^K\times \mathscr{C}_K ,
\end{equation}
where $\mathscr{C}_K := \{(\pi_1,\dots,\pi _K)\ |\ \pi_k \geq 0,\ \sum_{k=1}^K \pi_k =1\}$.
We define the probability distribution $\mathcal{P}_N$ by cutting off the $P_N$ in \eqref{log} as follows
\begin{align}
\log \mathcal{P}_N=\begin{cases}
\log P_N &\text{if } (\mu, \Lambda)\in \mathcal{PV}^K,\\
\infty&\text{otherwise}.
\end{cases}
\end{align}

In order to introduce the condition $\mathbf{B}$, we are required to set up some notations.
Let $\mathscr{Z}$ denotes the finite space
\[
\mathscr{Z}= \{1,\dots,K\}^N ,
\]
and for any element $\zeta=(\zeta_1,\dots,\zeta_N)\in \mathscr{Z},$ we associate a sequence 
\begin{align}\label{defz}
z:=\left\{z_{ik}\ |\ \ 1\leq i\leq N, \ 1\leq k\leq K\right\},
\end{align}
where 
\[
z_{ik}=\begin{cases}1&\text{if }\zeta_ i=k,\\ 
0&\text{otherwise}.\end{cases}
\]
We also set
\begin{align}\label{nkzikdef}
N_k=\sum_{i=1} ^N z_{ik},
\end{align}
and
\begin{equation}\label{defland}
\lambda(\zeta)=\min_{1\leq k\leq N} N_{k}.
\end{equation}
We consider a subset $\mathscr{Z}_0\subset \mathscr{Z}$ with the following properties
\begin{equation}
\mathscr{Z}_0=\left \{ \zeta \in \mathscr{Z}\ |\ \ \lambda(\zeta)\geq \lambda_0\right \},\label{land0}
\end{equation}
Let us define 
\begin{equation}\label{gmmdcm}
\mathcal{H}:=M \times \mathscr{Z}_{0}.
\end{equation}

\begin{condb}\label{condb1} The parameter $\lambda_0$ satisfies 
\[
\lambda_0 = O(N). 
\]
More precisely we assume that $\frac{\lambda_0}{N}> l_0$ for some positive constan $l_0$ which is independent of $N$. This means that we restrict the underlying parameter space to those partitions of 
the data into $K$ classes such that the minimum number of the data in each class is of $O(N)$. 
\end{condb}

Let $\bar{x}_k:= \frac{\sum_i z_{ik} x_i}{N_k}$ and $N_k=\sum_{n} z_{nk}$. To investigate the convexity of $-\log P_N$ in \eqref{log} we first need some simple calculation as below
\begin{align}\label{lam}
\sum _{n=1}^N z_{nk}\left(\mu _k - x_n \right)^T &\Lambda _k \left(\mu_k - x_n \right)\notag\\
=& \sum_{n=1}^N N_k \left(\mu _k - \frac{1}{N_k}\sum_{n=1}^N z_{nk} x_n \right)^T \Lambda _k \left(\mu_k - \frac{1}{N_k} \sum_{n=1}^N z_{nk} x_n\right)\notag \\
&+ \sum_{n=1}^N z_{n,k} x_n ^T \Lambda _k \left(z_{nk}x_n -\frac{1}{N_k} \sum_{n=1}^N z_{nk} x_n\right)\notag\\
=& \sum_{n=1}^N N_k \left(\mu _k - \frac{1}{N_k}\sum_{n=1}^N z_{nk} x_n \right)^T \Lambda_k \left(\mu_k - \frac{1}{N_k} \sum_{n=1}^N z_{nk} x_n\right) \notag\\
&+ \sum_{n=1}^N \left(z_{nk} x_n - \frac{1}{N_k} \sum_{n=1}^N z_{nk} x_n \right) ^T \Lambda _k \left(z_{nk} x_n -\frac{1}{N_k} \sum_{n=1}^N z_{nk} x_n\right)\notag\\
\quad =& \sum_{n=1}^N N_k \left(\mu_k - \bar{x}_k \right)^T \Lambda_k (\mu_k - \bar{x}_k)+ \sum_{n=1}^N z_{nk} (x_n- \bar{x}_{k})\Lambda_k (x_n - \bar{x}_k).
\end{align}

One can prove the following theorem,

\begin{theorem}\label{thcgmm}
Assume that conditions $\mathbf{A}$ and $\mathbf{B}$ are satisfied and the prior probability distributions $p(\Lambda_k)$ and $p(\mu_k)$ are of compact support and assume that $-\log p(\Lambda_k)$ and $-\log p(\mu_k)$ are $\mathfrak{C}-$convex for some constant $\mathfrak{C}\in \mathbb{R}$.
Then the potential $-\frac{1}{\lambda_0}\log \mathcal{P}_{N}$ is $C-$convex over each of the connectivity components of $\mathcal{H}$, where $C>0$ is a positive constant which does not depend on $N$.
\end{theorem}
\begin{proof}
According to \cite[Corollary 1]{bah}, we first note that given $X,Y \in S^P _{++}$ there exists an isometry $\mathcal{I}:S^P _{++} \rightarrow S^P _{++}$ such that $\mathcal{I}(X) = Id$ and $\mathcal{I} (Y) = \diag (e^r )$ where $r=(r_1,\dots, r_P) \in \mathbb{R} ^P$ and $\diag (e^r)$ denotes the diagonal matrix with diagonal coefficients equal to $(e^{r_1} ,\dots, e^{r_P})$. 
This isometry can be described through the action of an element $A$ of the group $GL(P)$ over $S^P _{++}$ defined by
\begin{equation}\label{act}
(Y, A) \rightarrow Y.A:= A^\dag Y A,
\end{equation}
where $Y \in S_{++} ^P$, $A\in GL(P)$ and $A^\dag$ is the transpose of $A$.\\
Let $\gamma$ be a geodesic in $S^P _{++}$ joining two elements $X,Y \in S^P_{++}$. Then one can find an orthogonal transformation 
$A\in O(P)$ such that $A^T X^{-1/2} \gamma X^{-1/2}A$ becomes a diagonalized path like \cite[Proposition 6]{bah}
\begin{equation}\label{diag}
t\rightarrow \diag(e^{r_1t},\dots,e^{r_P t}).
\end{equation}

If $\|r\|:=\sqrt{\sum_{i=1}^P r_i ^2} =1$, the geodesic will be of unit speed with respect to the geometry of $S^P_{++}$.
Thus, if we consider a linear change of coordinates on the data space $\mathbb{R}^P$ carried out both on $\mu_k$'s and $ x_i$'s, then it can be seen that the convexity of $\mathcal{P}_N$ is equivalent to the convexity along the paths in which the corresponding geodesic on $S^P _{++}$ is considered to be a diagonalized one.
For the sake of simplicity of notation, we now assume that the geodesic is diagonal. We assume that $\mu_k = (\mu_k ^1,\dots,\mu_k ^P)$ and
$\bar{x}_k= (\bar{x}_k ^1,\dots,\bar{x}_k ^P)$. For $k=1,\dots,K$ consider the unit speed geodesics 
\[
\gamma_k (t) = \bar{x}_k + a_k +t b_k, \quad 0\leq t\leq T^1 _k,
\]
in $\mu_k-$space with parameters $a_k=(a_k ^1,\dots,a_k ^P)$, and $b_k=(b_k ^1,\dots,b_k ^P)$ such that $\|b_k \|=1$.
Let also $\zeta_k$ be the unit speed geodesic
\[
\zeta_k (t) = \diag (e^{r_k ^1 t},\dots,e^{r_k ^P t}), \quad 0\leq t\leq T^2 _k, \quad 1\leq k \leq K.
\]
Then 
\begin{equation}\label{gdsg}
t\rightarrow \Gamma (t) :=\prod_{k=1} ^K \gamma_k (\alpha_k t) \times \zeta_k (\beta_k t) ,\quad\quad 0\leq t \leq T,
\end{equation}
where
\[
T=\sqrt{\sum_{k=1}^K (T_k ^1)^2 +\sum_{k=1}^K (T_k ^2)^2 }, 
\]
and for 
\begin{equation}\label{albet}
\alpha_k = \frac{T^1 _k}{T},\quad\text{ and }\quad
\beta_k =\frac{T^2_k }{T},\quad \text{}\ \ k=1,\dots,K,
\end{equation}
defines a geodesic in the space $\mathcal{PV}^K$.
Let the variable $u_k=(u_k ^1,{\dots},u_k ^P)$ be defined as $u^i _k:=\mu_k ^i- \bar{x}_k ^i,$ then by applying \eqref{lam} the restriction of $-\log\mathcal{P}_{N}$ to the geodesic \eqref{gdsg} is given by
\begin{align*}
-\log\mathcal{P}_{N}&=\sum_k N_k\left[\frac{1}{2} \sum_i \left(a^i _k +b^i _k \alpha_k t \right) ^2e^{r^i _k\beta_k t} +\frac{1}{2} \sum_{n,i} z_{nk} \left(x ^i _n-\bar{x} ^i _k\right)^2 e^{r^i _k \beta_k t} 
-\frac{1}{2} \sum_{n,i} z_{nk} r^i _k\beta_k t \right.\\
& +\left. \frac{1}{N_k}\!\sum_{i} \log p\left(a^i _k +b^i _k \alpha_k t+\bar{x}_k\right)+\frac{1}{N_k}\! \log p\left(\diag \left(e^{r_k\beta_k t}\right)\right) +\!\frac{1}{N_k}\log p\left(\pi_k\right)\! +\!const\right].
\end{align*}
Consequently,
\begin{align*}
-\frac{d^2}{dt^2}\log\mathcal{P}_{N}=&\sum_k N_k\left [\frac{1}{2} \sum_{i} \left [ 2\left(\alpha_k b_k ^i \right)^2 e^{r_k ^i \beta_kt} + 4\left(a_k ^i + b_k ^i \alpha_k t\right)b_k ^i r_k ^i \alpha_k \beta_k e^{r_k ^i\beta_k t}\right. \right.\\
&\left. +\left(a_k ^i + b_k ^i \alpha_k t\right)^2 \left(\beta_k r_k ^i \right)^2 e^{r_k ^i \beta_kt} \right]
+ \frac{1}{2}\sum_{n,i} z_{nk} \left(x_n ^i-\bar{x}_k ^i\right)^2 \left(r^i _k \beta_k \right)^2e^{r^i _k\beta_k t}\\
&\left. +\frac{1}{N_k}\frac{d^2}{dt^2} \left( \sum_{i} \log p\left(a^i _k +b^i _k\alpha_k t+\bar{x}_k\right)+ \log p\left(\diag \left(e^{r_k\beta_k t}\right)\right) \right) \right],
\end{align*}
and we have 
\begin{align*}
-\frac{d^2}{dt^2}\log\mathscr{P}_{N}=&\sum _k N_k \left[\left(\sum_i \left(b_k ^i \right)^2 e^{r_k ^i \beta_kt} \right)\alpha_k ^2
+ \left(2\sum_i \left(a^i _k + b^i _k \alpha_k t\right)b^i _k r^i _k e^{r^i _k \beta_k t}\right)\alpha_k \beta_k \right.\\
& +\left(\frac{t^2}{2}\sum_i \left(b_k ^i r^i _k \right)^2 e^{r^i _k \beta_k t} \right)\alpha_k ^2 \beta_k ^2
+\left(\sum_i a^i _k b^i _k t \left(r^i _k \right)^2 e^{r^i _k \beta_k t}\right)\alpha_k \beta_k ^2\\
& \left.+\frac{1}{2}\left(\sum_{n,i} z_{nk}\left (x_n ^i-\bar{x} ^i _k\right)^2 \left(r^i _k \right)^2e^{r^i _k\beta_k t}
+\sum_{i} \left(a^i _k r^i _k \right)^2e^{r^i _k\beta_k t} \right)\beta_k ^2\right. \\
& \left.+\frac{1}{N_k}\frac{d^2}{dt^2} \left( \sum_{i} \log p\left(a^i _k +b^i _k\alpha_k t+\bar{x}_k\right)+ \log p\left(\diag \left(e^{r_k\beta_k t}\right)\right) \right)\right] \\
\geq &\sum_k N_k \left(A_k \alpha_k ^2 + B_k \beta_k ^2+ C_k \alpha_k \beta_k\right) + \mathfrak{C},
\end{align*}
where
\[
\begin{split}
A_k &=\sum_i (b_k ^i )^2 e^{-\left| r_k ^i \right| t}, \\
\quad B_k& = \frac{1}{2}\sum_{n,i} z_{nk} \left(x_n ^i -\bar{x}_k ^i\right)^2 \left(r^i _k \right)^2 e^{ -\left|r^i _k\right| t}, \\
\quad C_k& = \sum_i 2 \left(a^i _k + b^i _k \alpha_k t\right)b^i _k r^i _k e^{r^i _k \beta_k t} -\left|\sum_i a^i _k b^i _k t \left(r^i _k \right)^2 e^{r^i _k \beta_k t}\right|,
\end{split}
\]
and $\mathfrak{C}$ introduced in Theorem \ref{thcgmm} is a constant such that
\[
\frac{d^2}{dt^2} \left(\sum_{k} \frac{1}{N_k}\sum_{i} \log p(a^i _k +b^i _k\alpha_k t+\bar{x}_k)+ \log p(diag (e^{r_k\beta_k t})) \right) \geq \mathfrak{C},
\]
from this and the fact that $t<2R$, 
\[
A_k \alpha_k ^2 + B_k \beta_k ^2+ C_k \alpha_k \beta_k\geq (\alpha_k ^2+\beta_k ^2)\nu,
\] 
for some positive $\nu$ which depends on $R$. To see this we note that given a symmetric matrix 
$A=\begin{pmatrix}a&c\\c&b \end{pmatrix},$
with $a,b>0,$ and $\det A>0,$ the smallest eigenvalue of $A$ is given by 
\[
\frac{2(ab-c^2)}{a+b + \sqrt{(a+b)^2-4(ab-c^2)}}\geq \frac{\det A}{\tr A},
\]
which means that
\[
|A|\geq \frac{\det A}{\tr A}.
\]
Also since $\lim_{b\rightarrow \infty} \frac{\det A}{ \tr A}= a,$ we see that for large enough $b$ we can lower estimate $|A|$ by $\frac{a}{2}$. From condition $\mathbf{A}$ we know that $A_k$ has a (positive) lower bound and $C_k$ has an upper bound both depending on $R$, and from condition $\mathbf{B}$ we know that $B_k$ grows to infinity as $N$ increases.

Since from \eqref{albet} we have $\sum_k \alpha_k ^2 + \beta_k ^2 =1$, therefore we can find a constant $C>0$ such that 
\[
\nu +\sum_k\frac{\mathfrak{C}}{N_k}>C>0.
\]
Hence we get to
\[
-\frac{d^2}{dt^2}\log\mathscr{P}_{N}\geq \min_k \{N_k\} \left[ \nu+\sum_k\frac{\mathfrak{C}}{N_k}\right].
\] 
This means that if we define $\Phi_N $ by $\Phi_N=-\frac{1}{\lambda_0}\log\mathcal{P}_N$ then we get
\[
-\log\mathcal{P}_N=\lambda_0 \Phi_N,
\] 
where $\lambda_0$ is defined by relation \eqref{defland} and we know that $\Phi_N$ is a $C-$convex map and $C$ can be chosen to be close to $\nu$ for large values of $N$.
\end{proof}

\begin{figure}[H]
\centering
\captionsetup{labelsep=none}
\begin{tikzpicture}[scale=0.5]
\def\verticalspacing{3}

\node at (0, .15\verticalspacing cm) {$\vdots$};

\foreach \i [count=\c] in {i,i',i''} {
\begin{scope}[yshift=-\c*\verticalspacing cm]
\draw[thick, black, domain=-1:1, samples=100] 
plot (\x, {4/(sqrt(2*pi)) * exp(-8*\x*\x/2)});
\filldraw[black] (0,0.142) circle (1pt);
\draw[thick, black, -] (-2.5,-1) -- (1.5,-1);
\draw[thick, black, -] (-2.5,-1) -- (-1,1);
\draw[thick, black, -] (1.5,-1) -- (2.5,1);

\node[below, font=\fontsize{8}{10}\selectfont] at (0.2,0.142) {$m_{\i}$};
\node[right, font=\fontsize{8}{10}\selectfont] at (3, 0.2) {$M \times \{ \i \}$};
\end{scope}
}

\begin{scope}[xshift=-3.5cm] 
\node[above, font=\fontsize{7}{10}\selectfont] at (-5.1,-5.9) {${Im}\,\mathfrak{M}$};
\draw[thick, black, domain=-6.2:-4.6, samples=100] 
plot (\x, {-6.5 + 0.8*exp(-(\x+6.5)^2/0.2) + 0.5*exp(-(\x+5.5)^2/0.4)});

\draw[thick, black, -] (-7.5,-7) -- (-4.5,-7);
\draw[thick, black, -] (-7.5,-7) -- (-6,-5);
\draw[thick, black, -] (-4.5,-7) -- (-3,-5);

\node[right, font=\fontsize{8}{10}\selectfont] at (-5,-7.5) {${M}$};

\draw[thick, black, ->] (-9.5,-6) -- (-7.18,-6);
\node[left, font=\fontsize{8}{10}\selectfont] at (-9.5,-6.175) { $\mathcal{A}$};

\node[above, font=\fontsize{8}{10}\selectfont] at (-8.35,-6) {$\mathfrak{\footnotesize M}$};

\draw[thick, black, ->] (-7.18,-6.35) -- (-9.5,-6.35);
\node[below, font=\fontsize{8}{10}\selectfont] at (-8.35,-6.4) {$\mathfrak{\footnotesize N}$};
\end{scope}

\node at (0, -3.5*\verticalspacing cm + 0.5) {$\vdots$};

\end{tikzpicture}
\caption{}
\label{fig:three_copies}
\end{figure}
%%%%%%%%%%%%%%

\section{Mean Field Variational Equations}\label{sec4}
The general setup for MFVBI applied to GMM consists of a hybrid discrete-continuous model in which the underlying space of the dataset has a product structure of the form $M\times \mathscr{Z}$ where $(M,g)$ is a complete Riemannian manifold and $\mathscr{Z}$ is a finite discrete set. The probability measure under investigation can be represented as $\frac{d\mu}{d\omega _g}=\frac{\exp \{-\lambda \Phi\}}{Z}$, where $\omega_g$ denotes the volume element associated with the Riemannian metric $g$. The restriction of the map $\Phi:M\times\mathscr{Z}\rightarrow \mathbb{R}$ to each slice $\Phi|_{M\times\{i\}}$ is assumed to be $C-$convex for all $i\in \mathscr{Z}$, where $C$ is a positive constant independent of $i$. Also $Z$ is the normalization constant such that $\mu$ belongs to the space $\mathcal{P} (M\times \mathscr{Z})$ of probability measures on $M\times \mathscr{Z}$.

The space of factorized probability measures denoted by $\mathscr{A}$ is defined to be 
\begin{equation}\label{defmathcalanew}
\mathscr{A}:=\mathcal{P}_2 (M)\times \mathcal{P}(\mathscr{Z}),
\end{equation}
where $\mathcal{P}_2(M)$ consists of the length space $(\mathcal{P}(M),W_2)$ of Borel probability measures $\mathcal{P}(M)$ on $M$ equipped with the 2-Wasserstein metric $W_2$. The variational
Bayesian problem consists of the following optimization problem
\begin{equation}\label{mfvbie}
\arg\min_{\nu\in \mathscr{A}}D_{KL} (\nu \|\mu).
\end{equation}

Since $\Phi|_{M\times \{i\}}$ is convex for each $i\in \mathscr{Z}$ there exists a unique point $(m_i,i)\in M\times \{i\}$ where the minimum of $\Phi|_{M\times \{i\}}$ occurs.

\begin{equation}\label{mipro}
m_i:=\argmin \Phi|_{M\times \{i\}}
\end{equation}
The absolute minimum of the application $\Phi$ is assumed to occur at $(m_{i_0}, i_0)$.

In our GMM problem the manifold $M$ is defined by relation \eqref{defm} and the discrete space $\mathscr{Z}$ 
equals $\mathscr{Z}_0$ defined by \eqref{land0}.
We assume that 
\begin{equation}\label{epla}
\lambda= \beta \lambda_0,
\end{equation}
where $\lambda_0$ is the same constant as given by \eqref{land0} which satisfies the condition $\mathbf{B}$.
The parameter $\beta$ can be viewed as a temperature parameter that will be fine-tuned later.

The application $\Phi$ is defined by 
\begin{equation}
\Phi:=\frac{1}{\lambda_0}\log \mathcal{P}_N.
\end{equation}
(See \eqref{defm} for defiition of $ \mathcal{P}_N$). 
Consider the probability measure $\mu_{}$ given by
\[
\frac{d\mu_{}}{d\omega_g}(\xi,i)=P_{}(\xi,i)=\frac{e^{-\lambda \Phi(\xi,i)}}{Z_{}},
\] 
where 

\[
Z_{}=\sum_{i\in \mathscr{Z}} \int_{M} e^{-\lambda \Phi (\xi,i)}d\omega_g
\]

is the normalization constant. Here $d\omega_g$ denotes the volume measure associated with the metric $g$. 
We note that the map $\Phi$ depends also on $\mathscr{Z}$ however by Theorem \ref{thcgmm} the positive constant $C$ representing its convexity coefficient, is independent of $\mathscr{Z}$.

If $\mu^1 _{} \times \mu^2 \in \mathscr{A}$ denotes the solution to the minimization problem \eqref{mfvbie} then we know that the following system of equations holds:

\begin{align}
\log \mu^2_{} (i) &=-\lambda\int_M \Phi (\xi,i) \frac{d\mu^1 _{}}{d\omega_g} d\omega_g-\log Z_{2},\label{mfvba1}
\\
\log \frac{d\mu^1 _{} }{d\omega_g}(\xi) &=-\lambda\sum_i \Phi (\xi,i) \mu^{2} _{} (i)-\log Z_{1}. \label{mfvba2}
\end{align}

\section{MFVBI and Marginal Partition Function}
For $i=1,\dots, |\mathscr{Z}|$ we define the parameters $u_i$ and $v_i$ as follows:
\begin{equation}\label{uivi}
u_i := \mu ^2 _{} (i) , \hspace{0.5cm}\text{ and }\hspace{0.5cm} v_i := \log u_i, 
\end{equation}
then from \eqref{mfvba2} we have 
\begin{equation}\label{eq38}
\frac{d\mu^1 _{} }{d\omega_g}(\xi) =\frac{e^{-\lambda\sum_i u_i \Phi (\xi,i)}}{Z_1}.
\end{equation}
Since $\sum_i u_i =1$, and the applications $\xi \rightarrow \Phi (\xi,i)$ is $C-$convex for all $i\in \mathscr{Z}_0$, the map 
\begin{align}\label{cconvexmu}
\xi \rightarrow \sum_i u_i \Phi (\xi,i),
\end{align}
will also become $C{-}$convex. Let $m\in M$ denote the minimum of the map $ \sum _i u_i \Phi (\xi,i)$:
\begin{align}\label{minmsigma}
m := \arg\min\limits_ {\xi\in M} \sum _i u_i \Phi (\xi,i),
\end{align}
By condition $\mathbf{B}$, the parameter $\lambda_0$ grows to infinity with $N$. Therefore 
the measure $\mu ^1 _{}$ concentrates at the single point $m$. More precisely $\mu^1 _{}$ will approach in probability towards the delta distribution $\delta (\xi-m)$.

On the other hand, by Laplace approximation, we know that for a convex function $\tilde{\Phi}$ and for any real continuous map $h:M\to \mathbb{R}$ we have
\begin{equation}\label{lapp}
\begin{split}
\int_{M} h(\xi)e^{-\lambda \tilde{\Phi}(\xi,i)} d\omega_g=& \left[\left|H_1(\tilde{\Phi}) (m_i,i)\right|^{-\frac{1}{2}}h(m_i) +O\!\left(\frac{1}{\sqrt{\lambda}}\right)\right]\left(\frac{2\pi}{\lambda}\right)^{\frac{d}{2}}e^{-\lambda\tilde{\Phi}(m_i, i)},
\end{split}
\end{equation}
where $H_1(\tilde{\Phi})$, denotes the Hessian of $\tilde{\Phi}$. Equivalently one can deduce that
\begin{equation}\label{lambdazlint}
\begin{split}
\int_{M} h(\xi)\frac{e^{-\lambda \tilde{\Phi}(\xi,i)}}{\tilde{Z}_L} d\omega_g=& h(m_i)+O\!\left(\frac{1}{\sqrt{\lambda}}\right) \left|H_1(\tilde{\Phi}) (m_i,i)\right|^{\frac{1}{2}},
\end{split}
\end{equation}
where $\tilde{Z}_L$ is defined by
$$
\tilde{Z}_L := \left|H_1(\tilde{\Phi}) (m_i,i)\right|^{-\frac{1}{2}} \left(\frac{2\pi}{\lambda}\right)^{\frac{d}{2}}e^{-\lambda\tilde{\Phi}(m_i, i)}.
$$

In particular we have
\begin{equation}\label{39}
\tilde{Z}_L=e^{-\lambda\tilde{\Phi}(m_i, i)+O(\log \lambda)}.
\end{equation}
From the above relation by setting $h\equiv 1$, it also follows that the normalization constant $\tilde{Z}:=\int e^{-\lambda \tilde{\Phi}(\xi,i)} $
is given by 
\begin{equation} \label{ztild}
\tilde{Z}:=\int e^{-\lambda \tilde{\Phi}(\xi,i)}=\tilde{Z}_L\left(1+ O\left(\frac{1}{\sqrt{\lambda}}\right)\right).
\end{equation}
So if we replace $\tilde{Z}_L$ in \eqref{lambdazlint} by $\tilde{Z}$ the right hand side of the relation \eqref{lambdazlint} can still be retained:
\begin{equation}\label{lambdazlintt}
\begin{split}
\int_{M} h(\xi)\frac{e^{-\lambda \tilde{\Phi}(\xi,i)}}{\tilde{Z}} d\omega_g=& h(m_i)+O\!\left(\frac{1}{\sqrt{\lambda}}\right) \left|H_1(\tilde{\Phi}) (m_i,i)\right|^{\frac{1}{2}}.
\end{split}
\end{equation}
Thus if we set 
\[
\tilde{\Phi}= \sum_i u_i \Phi (\xi,i),
\]
and
\[
v_i := \log \mu^2 (i),
\]
then from \eqref{mfvba1}, \eqref{minmsigma}, and the Laplace equation \eqref{lambdazlintt}, one can compute $v_i $
\begin{equation}\label{lapmu2}
\begin{split}
v_i &= \log \mu^2 _{ } (i)=\!-\lambda\!\int_M \Phi (\xi,i) \frac{d\mu^1 _{}}{d\omega_g} d\omega_g-\log Z_{2}=-\lambda\!\int_M \Phi (\xi,i) \frac{e^{-\lambda\sum_i u_i \Phi (\xi,i)}}{Z_1}d\omega_g-\log Z_{2}\\
& = -\lambda \Phi (m,i)+R_i-\log Z_2,
\end{split}
\end{equation}
where by \eqref{lambdazlintt}:
\begin{equation}\label{ri2}
R_i=O(\sqrt{\lambda}),
\end{equation}
and
\[
Z_2=\sum _i e^{-\lambda \Phi (m,i)+R_i}.
\]
Now by substituting \eqref{lapmu2} into \eqref{mfvba2} we obtain:
\begin{equation}\label{minmu}
\log \frac{d\mu ^1 _{}}{d\omega_g}= -\lambda \sum_i \Phi (\xi,i)e^{-\lambda \Phi (m,i)+ R_i-\log Z_2}-\log Z_1,
\end{equation}
where 
\begin{equation}\label{48}
Z_1= \int_M e^{-\lambda \sum \Phi (\xi,i)e^{-\lambda \Phi (m,i)+ R_i-\log Z_2}}d\omega_g.
\end{equation}

From the definition of $m$ in \eqref{minmsigma}, and the relations  \eqref{mfvba2} and \eqref{uivi}, the minimum of $-\log \frac{d\mu ^1 _{}}{d\omega_g}$ occurs at $\xi=m$. We also know from \eqref{cconvexmu} that $-\log \frac{d\mu ^1 _{}}{d\omega_g}$ is a convex function and therefore it admits a unique minimum which according to \eqref{minmu}, is a solution to the system of equations
\begin{equation}\label{dxi}
-\sum_i D_{\xi}\Phi (\xi,i)e^{-\lambda \Phi (m,i)+ R_i }=0.
\end{equation}
On the other hand any critical point of the map $\xi\rightarrow \sum_i e^{-\lambda \Phi (\xi,i)+ R_i }$ satisfies the equation
\[
-\sum_i D_{\xi}\Phi (\xi,i)e^{-\lambda \Phi (\xi,i)+ R_i }=0.
\]
Therefore the fact that $\xi = m$ satisfies the equation \eqref{dxi} is equivalent to say that $\xi =m$ is a critical point of the map
$\xi\rightarrow Z(\xi)$ defined by
\begin{equation}
Z(\xi):= \sum_i e^{-\lambda \Phi (\xi,i)+ R_i }.
\end{equation}

\begin{definition}
We call the map $\xi\rightarrow Z(\xi)$ defined by the above relation the marginal partition function of the GMM.
\end{definition}

We have thus proved the following theorem
\begin{theorem}\label{theoo2}
The maximum of $\log \frac{d\mu ^1 _{}}{d\omega_g}$ given by \eqref{minmu} with respect to $\xi$ is a critical point of the marginal partition function $\xi\rightarrow Z(\xi)$.
\end{theorem}

\section{Free Energy Associated with the Marginal Partition Function $Z(\xi)$}\label{zedxi}
The aim of this section is to effectively compute $Z(\xi)$ in probability with respect to the data $\{x_i\}_{i=1} ^N$ and for large values of $N$. We will show that there exists a map 
\[
\mathfrak{ M}: \mathcal{A}\rightarrow M,
\]
from the space of $K$ by $K$ Markov matrices $\mathcal{A}$ to $M$ such that in large $N$ limit we have
\begin{equation}\label{navas}
\frac{1}{\prod_{k=1} ^K (\tilde{N}_k)^{K-1}}\int (Z(\xi )) f(\xi) d\xi \rightarrow \int e^{-\lambda \big (\hat{\Phi} (A)+\psi (A)\big )+ O(\log \lambda)}\left(F(A)+O\left(\frac{1}{\lambda}\right)\right)d\mu_A,
\end{equation}
where $f:M\rightarrow \mathbb{R}$ is a test function and $F:\mathcal{A}\rightarrow \mathbb{R}$ is defined by $f= F\circ \mathfrak{M}^{-1}$. Also $\hat{\Phi}:\mathcal{A}\rightarrow \mathbb{R}$ given by
\begin{equation}\label{phihatag}
\hat{\Phi}(A):=\min_{\xi\in M}\hat{\phi}_A (\xi)=\hat{\phi}_A (\mathfrak{M} (A)),
\end{equation}
consists of a Legendre type transformation of $\hat{\phi}_A$ and $\xi\rightarrow \hat{\phi}_A (\xi)$ is derived from partitioning $\mathscr{Z}_0$ into some sub-classes in each of which $\Phi (\xi , i)$ concentrates around its average in that subclass (See \eqref{fiku}). 
The map $A\rightarrow \psi (A)$ asymptotically counts the number of elements in each sub-class and $d\mu_A$ denotes the Lebesgue measure on $\mathcal{A}$.

Based on the above relation \eqref{navas}, we define the effective partition function $Z_{eff}$ as follows
\begin{equation}\label{zeff0}
Z_{eff} (A) =e^{-\lambda \big (\hat{\Phi} (A)+\psi (A)\big )+ O(\log \lambda)}.
\end{equation}
In other words the relation \eqref{navas} shows that the marginal partition function $Z(\xi)$ converges in the weak sense towards a map concentrating around the subset $\mathfrak {M} (\mathcal{A})$. We will prove in Theorem \ref{the9} that the critical points of $Z(\xi)$ in the large $N$ limit converge towards the critical points of 
$A\rightarrow \big (\hat{\Phi} (A)+\psi (A)\big)$. This observation along with Theorem \ref{theoo2}.

The quantity $-\frac{1}{\lambda} \log Z_{eff}(A)$ is called according to physics literature the free energy of the system. 

\subsection{ Splitting the Data and U-Statistics}
In order to effectively compute the map $\xi\rightarrow Z(\xi)$ in large $N$ limit we observe that according to convexity and by Laplace approximation each of the terms $e^{-\lambda \Phi (\xi,i)+ R_i }$ behaves like a (non-normalized) delta distribution concentrated at the unique minimum point of $\xi\rightarrow \Phi (\xi ,i)$. Hence one can expect that the summation $\xi\rightarrow \sum_i e^{-\lambda \Phi (\xi,i)+ R_i }$ after normalization converges in the weak sense towards a map whose support lies on $\{m_i|i\in \mathscr{Z}\}$. In order to describe the asymptotic behavior of $Z(\xi)$ we first come up with a partition $\mathscr{Z}=\cup_{A_d\in \hat{\mathcal{A}}_d}\mathscr{Z}_{A_d}$ of $\mathscr{Z}$ into sub-classes $\mathscr{Z}_{A_d}$ in such a way that, in probability with respect to $\{x_i\}_{i=1} ^N$ 
the deviation of the maps $\xi\rightarrow \Phi (\xi , j)$ from their average for all $j\in \mathscr{Z}_{A_d}$ tends to zero with $N$.
The above assertion is proved by applying $U$-statistics Theorem \ref{chen}, and the limit of $\xi\rightarrow \Phi (\xi ,i))$ for $i\in \mathscr{Z}_{A_d}$ is denoted by 
$\xi\rightarrow \hat{\phi}_{A_d} (\xi) $.

Assume that the data $\{x_n\}_{n=1} ^N$ is sampled from the following mixture of normal distributions 
\begin{equation}\label{ptrue}
\tilde{p}_{N}(x)=\sum_{k=1} ^K \tilde{\pi}_k \ \mathcal{N}(x|\tilde{\mu}_k, \tilde{\Lambda}_k),
\end{equation}
which means that the true values of the parameters of our model consists of $\tilde{\pi}_k, \tilde{\mu}_k$ and $\tilde{\Lambda}_k$, $k=1,\dots,K$.
We also assume that the true classification of the data $\{x_n\}_{n=1} ^N$ is given by the parameters 
$\tilde{z}=\left\{\tilde{z}_{ik} | \ 1\leq i \leq N, \ 1\leq k \leq K\right\}.$
Let us define
\[
A_{z}(k,k' ):= \left|\left\{i\ |\ z_{ik}=\tilde{z}_{ik'}=1 \right\}\right|.
\] 
In other words, $A_{z}( k, k')$ enumerates the number of the elements of the data which belong to the (true) $k'$-th class while they are classified as being in the $k$-th class according to the partition induced by $z$ (see the definition \ref{defz}). Consequently, we have 
\begin{equation}\label{mata1}
\sum _{k'} A_{z}( k, k')= N_k,
\end{equation}
and
\begin{equation}\label{mata}
\sum _{k} A_{z}( k, k')=\tilde{ N}_{k'},
\end{equation}
where $\tilde{N}_k'=\sum_i \tilde{z}_{ik'}$ represents the true number of the data in $k'$-th class and $N_k=\sum_i z_{ik}$ denotes the number of the data in $k$-th class
according to the hypothetical classification given by $z$.
Let $\hat{\mathcal{A}}_d$ denote the collection of $K\times K$ matrices with non-negative integer coefficients satisfying the relation \eqref{mata}:
\begin{equation}\label{Addef}
\hat{\mathcal{A}}_d\!:=\!\left\{\!\hat{A}_d=\left[\hat{A}_{d}(k,k')\right]_{K\times K}\in M_{K}\! (\mathbb{N}\cup \{0\}) \left| \sum_{k=1}^K \hat{A}_{d}(k,k') =\tilde{N}_{k'}, \text{ for } 1\leq k\leq K\right.\right\}.
\end{equation}
Associated with each $\hat{A}_d \in \hat{\mathcal{A}}_d$ we define a subset $\mathscr{Z}_{\hat{A}_d}\subset\mathscr{Z}$ as follows:
\[
\mathscr{Z}_{\hat{A}_d}=\left\{z\ |\ \ A_{z}(k,k')=\hat{A}_{d}( k,k') \text{ for all } 1\leq k,k'\leq K\right\}.
\]
Hence $\mathscr{Z}=\cup_{\hat{A}_d\in \hat{\mathcal{A}}_d} \mathscr{Z}_{\hat{A}_d}$ provides a partition of $\mathscr{Z}$ into a union of disjoint subclasses $\mathscr{Z}_{\hat{A}_d}$ for $\hat{A}_d\in \hat{\mathcal{A}}_d$. Thereof we can split up the summation $\sum_{i\in \mathscr{Z}} \Phi (\xi, i)$ as follows:
\begin{equation}
\sum_{i\in \mathscr{Z}}e^{-\lambda \Phi (\xi , i)+R_i}= \sum_{\hat{A}_d \in \hat{\mathcal{A}}_d} \sum_{i\in \mathscr{Z}_{\hat{A}_d}}e^{-\lambda \Phi (\xi , i)+R_i}.
\end{equation}
The fundamental feature of the subclasses $\mathscr{Z}_{\hat{A}_d}$ is that by U-statistics Theorem \ref{chen} for a fixed $\hat{A}_d $ the maps $\{ \xi\rightarrow \Phi(\xi,z)\}_{z\in \mathscr{Z}_{\hat{A}_d}}$ are close together with respect to compact open topology and approach their average in probability with respect to the choice of the data $\{x_i\}_{i=1} ^N$. 

More precisely if we define the rational numbers $\alpha_{N}(k,k')\in \mathbb{Q}\cap [0,1]$ by 
\begin{equation}\label{alfa1}
\alpha_{N} (k,k'):=\frac{\hat{A}_{d}(k,k')}{\tilde{N}_{k'}}.
\end{equation}

Then, according to Proposition \ref{prop6} in the Appendix \ref{app4}, the average of $\lambda\Phi|_{M\times \mathscr{Z}_{\hat{A}_d}}$ with respect to $\tilde{p}_N$ (defined by \eqref{ptrue}) equals
\begin{align} 
\lambda_0 \left\langle \Phi |_{M\times \mathscr{Z}_{\hat{A}_d}} \right\rangle_{\tilde{p}_{N} }=&-\frac{1}{2} \sum_{k=1}^K \left(\sum_{k'=1}^K \tilde{N}_{k'}\alpha_{N} (k,k')\tilde{\pi}_{k'}\sum_{1\leq i,j\leq P} \tilde{\Lambda}^{ij} _{k'} \Lambda_{k,ij}\right.\notag\\
&\quad \quad \quad \quad \quad \quad \left.+\sum_{k'=1}^K \tilde{N}_{k'} \alpha_{N} (k,k')\tilde{\pi}_{k'} (\mu_k -\tilde{\mu}_{k'})^T\Lambda_{k} (\mu_k -\tilde{\mu}_{k'}) \right)\notag \\
\quad &+ \left [ \sum_{k=1}^K \sum_{k'=1}^K \tilde{N}_{k'} \tilde{\pi}_{k'}\alpha_N (k,k')\left( \log \pi _k
+ \frac{1}{2} \log |\Lambda _k | \right) \right ]+ \sum _{k=1} ^K \log p (\mu _k )\notag\\
&+ \sum _{k=1} ^K \log p (\Lambda _k ) +\log p (\pi) +\tilde{C}.
\label{fiku2}
\end{align}

It follows from \eqref{mata} that $[\alpha_{N} (k,k')]_{1\leq k,k'\leq N}$ is in fact a Markov matrix.
We consider the space of all Markov $K\times K$ matrices 
\begin{equation}\label{mathcala}
\mathcal{A}\!=\!\left\{\! A=\!\left[\alpha (k,k')\right]\!\in M_{K\times K} (\mathbb{R}) \bigg|\sum_{k =1} ^K \alpha (k,k')=1, \ \alpha (k,k')\geq 0 \text{ for } 1\leq k, k' \leq K \right\}.
\end{equation}
If we set
\begin{equation}\label{iws1}
\mathcal{A}_d:=\{A=[\frac{\hat{A}_d (k,k')}{\tilde{N}_{k'}}]_{k,k'}\in M_{K\times K} (\mathbb{Q})| [\hat{A}_d (k,k')]_{k,k'}\in \hat{\mathcal {A}_d}\}
\end{equation}
then $\mathcal{A}_d$ is a lattice in $\mathcal{A}$ and the map
\begin{equation}\label{j1}
\hat{} :\mathcal{A}_d\rightarrow \hat{\mathcal{A}}_d,\hspace{0.5cm} (A_d)\rightarrow \hat{A}_d
\end{equation}
\begin{equation}\label{j2}
\hat{A}_d=[\hat{A}_d (k,k')]_{k,k'} \hspace{0.5cm} A_d=[\frac{\hat{A}_d (k,k')}{\tilde{N}_{k'}}]_{k,k'},
\end{equation}
establishes a one to one correspondence between $\hat{\mathcal{A}}_d$ and $ \mathcal{A}_d$.
Also based on \eqref{fiku2} for any $A\in \mathcal{A}$ and $1\leq k \leq K$ we define $\hat{\phi}_{A} (\xi)$ by
\begin{equation}\label{fiku}
\begin{split}
-\lambda_0\hat{\phi}_{A} (\xi):=&-\frac{1}{2} \sum_{k=1}^K \left(\sum_{k'=1}^K \tilde{N}_{k'}\alpha (k,k')\sum_{1\leq i,j\leq P} \tilde{\Lambda}^{ij} _{k'} \Lambda_{k,ij}\right.\\
&\quad\quad\quad\quad\quad\quad\left.+\sum_{k'=1}^K \tilde{N}_{k'}\alpha (k,k') (\mu_k -\tilde{\mu}_{k'})^T\Lambda_{k} (\mu_k -\tilde{\mu}_{k'}) \right) \\
\quad &+ \left [ \sum_{k=1}^K \sum_{k'=1}^K \tilde{N}_{k'} \alpha_N (k,k')\left( \log \pi _k
+ \frac{1}{2} \log |\Lambda _k | \right) \right ]\\
&+ \sum _{k=1} ^K \log p (\mu _k )+ \sum _{k=1} ^K \log p (\Lambda _k ) +\log p (\pi) +\tilde{C},
\end{split}
\end{equation}
where $\lambda_0$ was introduced in \eqref{defland}. Also, according to \eqref{fiku2} and \eqref{fiku} and from the U-statistics theory \cite{Chen} one can see that for $i\in \mathscr{Z}_{A_d}$
\begin{equation}\label{fdl}
\Phi (\xi ,i ) = \hat{\phi}_{A_d} (\xi)+ O (\frac{1}{\sqrt{\lambda}}).
\end{equation}
Let also 
\begin{equation}
\mathfrak{ M}: \mathcal{A}\rightarrow M,
\end{equation}
be the application which assigns to $A\in \mathcal{A}$ the point 
\[
\mathfrak{ M}(A):=\left(\Lambda_1 (A), \dots, \Lambda_K (A), \mu_1 (A), \dots, \mu_K (A), \pi_1 (A), \dots, \pi_K (A)\right)\in M,
\]
at which the minimum of $\hat{\phi}_{A} (\xi)$ occurs. Note that by Theorem \ref{thcgmm} we know that this minimum is unique. 
Next, we define the map $\hat{\Phi}:\mathcal{A}\rightarrow \mathbb{R}$, by the following Legendre type transformation:
\begin{equation}\label{phihatag}
\hat{\Phi}(A):=\min_{\xi\in M}\hat{\phi}_A (\xi)=\hat{\phi}_A (\mathfrak{M} (A)).
\end{equation}
We remark that, as is standard for Legendre transformation since $\hat{\phi}_A(\xi)$ is linear with respect to $A$, the map $A\rightarrow \hat{\Phi}(A)$ will be concave.

\subsection{Estimation of the Number of the Elements in Each Subclass $\mathscr{Z}_{\hat{A}_d}$}
For simplicity of notation we replace $\alpha_{N}(k,k')$ by $\alpha_{k,k'}$.
By Stirling's approximation 
\[
N! \sim \left(\frac{N}{e}\right)^N \sqrt{2\pi N},
\]
thus, we have:
\begin{align}
|\mathscr{Z}_{\hat{A}_d}|&=\prod\limits_{k'=1}^K {\tilde{N}_{k'}\choose \hat{A}_d (1,k'), \dots ,\hat{A}_d (K,k')}=
\prod_{k'=1}^K\frac{\tilde{N}_{k'}!}{\hat{A}_d (1,k')!\dots \hat{A}_d (K,k')!}=
\prod_{k'=1}^K\frac{\tilde{N}_{k'}!}{\prod\limits_{k=1}^K \hat{A}_d (k,k')!}\notag\\
&\simeq
\prod_{k'=1}^K\frac{\left(\frac{\tilde{N}_{k'}}{e}\right)^{\tilde{N}_{k'}}\sqrt{2\pi \tilde{N}_{k'}}}{\prod\limits_{k=1}^K 
\left(\frac{\hat{A}_d (k,k')}{e}\right)^{\hat{A}_d (k,k')}\sqrt{2\pi\hat{ A}_d (k,k')}}
=
\prod\limits_{k'=1}^K\frac{\left(\tilde{N}_{k'}\right)^{\tilde{N}_{k'}}\sqrt{2\pi \tilde{N}_{k'}}}{\prod\limits_{k=1}^K 
\left(\hat{A}_d (k,k')\right)^{\hat{A}_d (k,k')}\sqrt{2\pi \hat{A}_d (k,k')}} \notag\\
&=
(2\pi) ^{\frac{K-K^2}{2}}\prod\limits_{k'=1}^K \frac{{\tilde{N}_{k'}}^{\tilde{N}_{k'}}}{\prod \limits_{k=1}^K (\hat{A}_d (k,k'))^{\hat{A}_d (k,k')}} \times \prod_{k'=1}^K\frac{{\tilde{N}_{k'}}^{\frac{1}{2}}}{\prod\limits _{k=1}^K (\hat{A}_d (k,k'))^{\frac{1}{2}}}\notag\\
&= (2\pi) ^{\frac{K-K^2}{2}}\prod\limits_{k'=1}^K\left( \prod\limits_{k=1}^K (\alpha_{k,k'})^{\alpha _{k,k'}}\right)^{-\tilde{N}_{k'}} \times \prod\limits _{k'=1}^K\left( (\tilde{N}_{k'})^{\frac{1-K}{2}} \prod\limits_{k=1}^K (\alpha_{k,k'})^{\frac{-1}{2}}\right). \label{iternumz}
\end{align}
We define $\psi$ as follows:
\begin{equation}\label{sydef}
-\lambda \psi=\log \left| \mathscr{Z}_{\hat{A}_d}\right|=-\sum_{k,k'} \tilde{N}_{k'}\alpha_{k,k'}\log \alpha_{k,k'}+O(\log \tilde{N}_{k'}),
\end{equation}
so 
$\psi$ is a smooth map over $\mathcal{A}$. 

\subsection{Computation of $Z_{eff} (A)$ in the Weak Sense}
If we assume that $\mathfrak{M}:\mathcal{A}\rightarrow M$ is injective then for any map $f: \mathfrak{M}(\mathcal{A})\rightarrow \mathbb{R}$ there exists $F:\mathcal{A}\rightarrow \mathbb{R}$ such that $f= F\circ \mathfrak{M}^{-1}$. In order to describe the asymptotic behavior of the marginal partition function $Z(\xi)$ we need to carry out a proper normalization. Here in this section, we do a primary computation to obtain an intuition of what we expect of the limit $\lim_{\lambda\rightarrow \infty}Z(\xi)$ in the weak sense. A more rigorous treatment of this limit will be performed elsewhere. We define
\[
\tilde{Z}_0:=\sum_{A_d \in \mathcal{A}_d} e^{-\lambda (\hat{\phi}_A (\mathfrak{M} (A_d) +\psi (A_d))+O(\sqrt{\lambda})},
\]
and we can study the normalized partition function $\frac{Z}{\tilde{Z}_0}$. But here for the sake of simplicity, we divide $Z$ by $2^{D/2}\prod_{k=1} ^K (\tilde{N}_k)^{K-1}$ which consists of the volume of each of the cells of the lattice $\mathcal{A}_d$ inside $\mathcal{A}$. We have also set $D=K(K-1)$.

By Laplace approximation in \eqref{lapp}, the number $|\mathscr{Z}_{A_d}|$ of iterations approximated in \eqref{iternumz} as well as the relation \eqref{fdl}, we have
%%%%%%%%%%%%%%
\begin{align}\label{calint1}
&\frac{1}{2^{D/2}\prod_{k=1} ^K (\tilde{N}_k)^{K-1}}\int Z(\xi ) f(\xi) d\xi\notag\\
&\quad \quad=\frac{1}{2^{D/2}\prod_{k=1} ^K (\tilde{N}_k)^{K-1}}\int \big ( \sum_{i\in \mathscr{Z}} e^{-\lambda \Phi(\xi ,i)+R_i} \big)f(\xi)d\xi\notag\\
%\quad & = \int \big (\sum_{i\in \mathscr{Z}} e^{-\lambda \Phi(\xi ,i)+R_i}\big ) f'(\xi) d\xi \\
&\quad \quad = \frac{1}{2^{D/2}\prod_{k=1} ^K (\tilde{N}_k)^{K-1}}\int ( \sum_{\hat{A}_d \in \hat{\mathcal{A}}_d} \sum_{i\in \mathscr{Z}_{\hat{A}_d}}e^{-\lambda \Phi (\xi , i)+R_i}) f(\xi) \big )d\xi\notag\\
&\quad \quad = \frac{1}{2^{D/2}\prod_{k=1} ^K (\tilde{N}_k)^{K-1}}\sum_{\hat{A}_d \in \hat{\mathcal{A}}_d} \sum_{i\in \mathscr{Z}_{\hat{A}_d}} \big (f(m_i)+O(\frac{1}{\sqrt{\lambda}})\big )e^{-\lambda \Phi(m_i, i)+ O(\log \lambda)+R_i}
\notag\\
&\quad \quad \rightarrow \frac{1}{2^{D/2}\prod_{k=1} ^K (\tilde{N}_k)^{K-1}}\sum_{\hat{A}_d\in \hat{\mathcal{A}}_d} | \mathscr{Z}_{\hat{A}_d}|e^{-\lambda\hat{\phi} _{A_d}(\mathfrak{M} (A_d))+ O(\sqrt{\lambda})+R_i} \big ( f(\mathfrak{M}(A_d)+O(\frac{1}{\sqrt{\lambda}})\big )\notag\\
&\quad \quad = \frac{1}{2^{D/2}\prod_{k=1} ^K (\tilde{N}_k)^{K-1}}\sum_{\hat{A}_d\in \hat{\mathcal{A}}_d} e^{-\lambda\big (\hat{\phi}_{A_d}(\mathfrak{M} (A_d)) +\psi (A_d)\big )+ O(\sqrt{\lambda})+R_i} \big (F(A_d)+O(\frac{1}{\sqrt{\lambda}})\big )\notag\\
&\quad \quad \rightarrow \int e^{-\lambda \big (\hat{\Phi} (A)+\psi (A)\big )+ O(\sqrt{ \lambda})+R_i}(F(A)+O(\frac{1}{\sqrt{\lambda}}))d\mu_A,
\end{align}
where $d\mu_A$ is the Lebesgue measure induced on $\mathcal{A}$ as a subspace of $\mathbb{R}^{K^2},$ and where in the third line we are using U-statistics Theorem \ref{chen} in the large $\lambda$ limit through the relation \eqref{fdl}. The limit in the last line is incorrect in the form it is written. We prove in a separate work that after appropriate normalization
\[
\lim_{\lambda\rightarrow \infty} \int \frac{Z(\xi)}{\tilde{Z}_0 } f(\xi) d\xi=\sum c_i f(\mathfrak{M}(A_i))
\]
where $A_1, \dots, A_m$ are the points where the absolute minimum of $\hat{\Phi} (A)+\psi (A)$ occurs and $c_1, \dots, c_m$ are some constants which only depend on $\hat{\Phi} (A)+\psi (A)$. Here $\tilde{Z}_0 $ is the normalization constant introduced earlier.

Therefore, we consider $\hat{\Phi} (A)+\psi (A)$ as free energy. 
Based on the above equation we define the effective partition function $Z_{eff} (A)$ as follows
\begin{equation}\label{zeff}
Z_{eff} (A) =e^{-\lambda \big (\hat{\Phi} (A)+\psi (A)\big )}.
\end{equation}

\section{Main Theorems and a Reduced MFVI Equation}\label{MainT}

The following theorem demonstrates that the connection established in Theorem \ref{theoo2} between the mode of the solution to MFVBI and critical points of partition function 
is retained equally for the effective partition function. It also characterizes the asymptotic behavior of $\frac{d\mu^1}{d\omega_g}$ in large $\lambda$ limit.

\begin{theorem}\label{the9} If the maximum of the right-hand side of \eqref{minmu}
\begin{equation}\label{minmu22}
\log \frac{d\mu ^1 _{}}{d\omega_g}= -\lambda \sum_i \Phi (\xi,i)e^{-\lambda \Phi (m,i)+ R_i-\log Z_2}-\log Z_1,
\end{equation}

with respect to $\xi$ occurs at a point which corresponds to an interior point of $\mathcal {A}$ under the map $\mathfrak{M}$, then in large $\lambda$ limit this point converges towards a critical point of $A\rightarrow (\hat{\phi}_A (\mathfrak{M} (A) +\psi (A)) $. If the maximum occurs at the interior point of any of sub-complexes of $\mathcal{A}$ it will coincide with a critical point of the restriction of $A\rightarrow (\hat{\phi}_A (\mathfrak{M} (A) +\psi (A)) $ to that sub-complex. 
\end{theorem}
\begin{proof}
As before by applying \eqref{fdl}, one can see that
\begin{equation}\label{1233}
\sum_i \Phi (\xi,i)e^{-\lambda \Phi (m,i)+ R_i-\log Z_2}\simeq \sum_{A_d}\frac{e^{-\lambda (\hat{\phi}_{A_d} (m) +\psi (A_d))+O(\sqrt{\lambda})} \hat{\phi}_{A_d} (\xi) }{\tilde{Z}_1}
\end{equation}
\begin{equation}\label{ztild1}
\tilde{Z}_1=\sum_{A_d} e^{-\lambda (\hat{\phi}_{A_d} (m) +\psi (A_d))+O(\sqrt{\lambda})} 
\end{equation}
where $m$ is defined in \eqref{lambdazlintt}. 
In fact from \eqref{mfvba1} we have 
\[
Z_2= \sum _i e^{-\lambda \int _M \Phi (\xi, i) \frac{d\mu ^1}{d\omega_g} d\omega_g}.
\]
Thus, by \eqref{eq38} and \eqref{minmsigma} and from Laplace approximation one can see that
\[
Z_2\simeq \sum_i e^{-\lambda \Phi (m,i)+ O(\sqrt{\lambda}) } .
\]
This combined with \eqref{fdl} justifies the introduction of $\tilde{Z}_1$ in \eqref{1233} and \eqref{ztild1}.

We can assume that $ (\hat{\phi}_A (m) +\psi (A))+\frac{1}{\lambda}O(\sqrt{\lambda})\geq 0$ and $0$ occurs as its minimum.
Also we assume that $\hat{\phi}_A (\xi)\geq 1 $. Both of these assumptions can be established by a constant shift. Let $B$ be the point where the minimum of $\hat{\phi}_{A} (m) +\psi (A))+\frac{1}{\lambda}O(\sqrt{\lambda})$ occurs:
\begin{equation}\label{124}
B:=\argmin_ A \big ( \hat{\phi}_{A} (m) +\psi (A))+\frac{1}{\lambda}O(\sqrt{\lambda})\big )
\end{equation}
We determine $\epsilon_0$ in such a way that for $A$ satisfying $\|A- B \| \geq \epsilon_0$ and for all $A_{d_0}$ which have minimum distance with respect to $B$ we have 
\begin{equation}\label{137}
(\hat{\phi}_A (m) +\psi (A))+\frac{1}{\lambda}O(\sqrt{\lambda})\geq \delta+ (\hat{\phi}_{A_{d_0}} (m) +\psi (A_{d_0}))+\frac{1}{\lambda}O(\sqrt{\lambda}).
\end{equation}
Here $\delta$ is an arbitrary positive constant. 
We have $\tilde{Z}_1 \geq e^{-\lambda (\hat{\phi}_{A_{d_0}} (m) +\psi (A_{d_0}))+O(\sqrt{\lambda})} $ where $A_{d_0}$ is of least distance to $B$. So by applying \eqref{137} one can deduce that for $A_d$ satisfying $\|A_d-B\|\geq \epsilon_0$ and $\|A_d-A_{d_0}\|\geq \epsilon_0$ for all $A_{d_0}$ having the least distance with respect to $B$
\begin{align*}
\frac{e^{-\lambda (\hat{\phi}_{A_d} (m) +\psi (A_d))+O(\sqrt{\lambda})} \hat{\phi}_{A_d} (\xi) }{\tilde{Z}_1}\leq e^{-\lambda \delta}\max_{A} \hat{\phi}_A (\xi) .
\end{align*}

Therefore $\lim_{\lambda\rightarrow \infty} \sum_{\{A_d| \|A_d-B\|\geq \epsilon_0\}} \frac{1}{\tilde{Z}_0} e^{-\lambda (\hat{\phi}_{A_d} (m) +\psi (A_d))+O(\sqrt{\lambda})} \hat{\phi}_{A_d} (\xi)=0$. (Since the cardinality of the set $\{A_d| \|A_d-B\|\geq \epsilon_0\}$ has a polynomial growth with respect to $\lambda$.) So 
\begin{align}\label{1277}
\sum_{A_d} &\frac{1}{\tilde{Z}_0} e^{-\lambda (\hat{\phi}_{A_d} (m) +\psi (A_d))+O(\sqrt{\lambda})} \hat{\phi}_{A_d} (\xi) \simeq \notag \\
&\sum_{\{A_d| \|A_d-B\|\leq \tilde{\epsilon}_0\}} \frac{1}{\tilde{Z}_0} e^{-\lambda (\hat{\phi}_{A_d} (m) +\psi (A_d))+O(\sqrt{\lambda})} \hat{\phi}_{A_d} (\xi) \simeq \hat{\phi}_{B} (\xi) \sum_{\{A_d| \|A_d-B\|\leq \tilde{\epsilon}_0\}} \frac{1}{\tilde{Z}_0} e^{-\lambda (\hat{\phi}_{A_d} (m) +\psi (A_d))+O(\sqrt{\lambda})}
\end{align}
Here $\tilde{\epsilon}_0=\epsilon_0+d(B, \mathcal{A}_d)$.
Therefore from \eqref{1233} and \eqref{1277} it follows that the minimum 
\[
\argmin_{\xi} \sum_i \Phi (\xi,i)e^{-\lambda \Phi (m,i)+ R_i-\log Z_2},
\]
approaches $\xi=\mathfrak{M}(B)$ as $\lambda$ grows. In other words 
\begin{equation}\label{126}
m\simeq \mathfrak{M}(B).
\end{equation}
Also by the definition of $B$ in \eqref{124}
\begin{equation}\label{119}
\hat{\phi}_A (\mathfrak{M}(B))+\hat{\psi} (A)+\frac{1}{\lambda}O(\sqrt{\lambda})\geq \hat{\phi}_B (\mathfrak{M}(B))+\hat{\psi} (B)+\frac{1}{\lambda}O(\sqrt{\lambda}), \hspace{1cm}\text{ for all } A.
\end{equation}
Due to relations \eqref{124} and \eqref{126} one can also see that, 
\begin{equation}\label{130}
D_{A}\big ( \hat{\phi}_A (\xi)+\psi (A) \big ) |_{\xi =\mathfrak{M}(B), A= B} \simeq 0
\end{equation}
From the definition of $\mathfrak{M}$ we also know that
\begin{equation}\label{131}
D_{\xi} \hat{\phi}_B (\xi)|_{\xi=\mathfrak{M}(B)}\simeq 0.
\end{equation}
Now by \eqref{130} and \eqref{131} by applying chain rule one can conclude that $B$ is a critical point of $A\rightarrow \hat{\phi}_A (\mathfrak{M}(A))+\hat{\psi}(A)$
\begin{equation}
D_A \big ( \hat{\phi}_A (\mathfrak{M}(A))+\psi (A) \big )|_{A=B} \simeq 0.
\end{equation}
\end{proof}
The computation leading to \eqref{1277} does not depend on the precise values of the coefficients $\Phi (\xi,i)$
in the equation \eqref{minmu22}. Thus, from this computation one can deduce that the parobability measure $\mu_1$ satisfying \eqref{mfvba2} 
concentrates in large $\lambda$ limit at those $i$ which belonge to sub classes $\mathscr{Z}_{A_d}$ with $A_d$ approching $B$.
Here $B$ saisfies the equation \eqref{124}.
Also by \eqref{126} we have $\mathfrak{M} (B)=m$. If we define $\mathfrak{N}: M\rightarrow \mathcal{A}$ as the Legendre transform of $\xi\rightarrow \hat{\phi}_A (\xi ) +\hat{\psi} (A) $
\begin{equation}
\mathfrak{N} (\xi ):=\argmin_A \big ( \hat{\phi}_A (\xi ) +\hat{\psi} (A) \big).
\end{equation}
From \eqref{126} and \eqref{124} we get 
\begin{equation}
\mathfrak{N} (\mathfrak{M}(B))\simeq B.
\end{equation}
In other words, the probability measures $\mu^1$ and $\mu^2$, can be approximated in large $\lambda$ limit by the following measures
\begin{equation}\label{dirac}
\frac{d\mu^1}{d\omega_g}\simeq \delta (\xi-m), \hspace{0.5cm} \mu^2 (i) \simeq \frac{1}{|\mathscr{Z}_{\hat{A}_{d_1}}|} \sum_{j\in \mathscr{Z}_{\hat{A}_{d_1}}} \delta(i-j),
\end{equation}
where $A_{d_1}$ is one of the points of the lattice $\mathcal{A}_d$ possessing the minimum distance with respect to $B$.
\[
A_{d_1}\simeq B,
\]
In fact, if we define 
\[
\mathfrak{J}:\mathscr{Z}\rightarrow \mathcal{A},
\]
by
\[
\mathfrak{J}(i)= A_d, \hspace{1cm} \text{for } i\in \mathscr{Z}_{\hat{A}_d}.
\]
(See \eqref{j1} for the definition of $\hat{}$ ) then we have 
\[
\mu^2\simeq \mathfrak{J}^* \left [\delta (A-A_{d_1}) \right ].
\]
\begin{theorem}\label{newmfvba}
The two Dirac measures given by the equation \eqref{dirac} constitute a solution to the MFVBA equations \eqref{mfvba1} and \eqref{mfvba2} 
iff $B$ and $m$ satify the following system of equations
\begin{equation} \label{bonyadi}
\mathfrak{N} (\mathfrak{M}(B))\simeq B, \hspace{0.5cm} \mathfrak{M} (B)=m
\end{equation}
\end{theorem}
\begin{proof}
The proof of "only if" direction follows from the discussion provided above. For the "if" direction if we substitute $\sum_{i\in \mathscr{Z}_{A_{d_1}}}\delta (i-j)$ into the equation \eqref{mfvba2}, then one gets
\[
\frac{d\mu^1}{d\omega_g} \simeq\frac{ e^{-\lambda \sum_{i\in \mathscr{Z}_{A_{d_1}}} \Phi (\xi,i)}}{Z_1}.
\] 
Hence by the choice of $ \mathscr{Z}_{A_{d_1}}$ for $i\in \mathscr{Z}_{A_{d_1}}$ we have $m_i\simeq \mathfrak{M}(B)$. We also note that from the definition of $m$ in \eqref{lambdazlintt} we have 
\begin{equation}
m\simeq m_i\simeq \mathfrak{M}(B).
\end{equation}
Therefore, from \eqref{mfvba1} we get
\begin{equation}\label{minmu22}
\log \frac{d\mu ^1 _{}}{d\omega_g}= -\lambda \sum_i \Phi (\xi,i)e^{-\lambda \Phi (m,i)+ R_i-\log Z_2}-\log Z_1.
\end{equation}
Then by the same argument as in the proof of Theorem \ref{the9}, $\frac{d\mu^1}{d\omega_g}$ is concentrated at 
$\mathfrak{M} (B')$ where $B'$ is given by
\begin{equation}
B'=\mathfrak{N} (m): =\argmin_ A \big ( \hat{\phi}_{A} (m) +\psi (A))+\frac{1}{\lambda}O(\sqrt{\lambda})\big )
\end{equation}

But from \eqref{bonyadi} we have $\mathfrak{N} (m)=B$, hence $B=B'$ and this shows the compatibility required for 
the system of equations \eqref{mfvba1} and \eqref{mfvba2}. 
\end{proof}

From the relation \eqref{fiku2} it is not difficult to observe that the minimum of $\hat{\phi}_A (\xi)$ coincide with the right parameters 
iff hte matrix $A$ is the identity matrix $I_k$ or one of the matrices obtained by an elementary row operation on $I_K$.
In fact by starting the minimization with respect to $\mu_k$'s it can be seen that the only values for $A$ yielding $\mu_k=\tilde{\mu_k}$ (upto a permutation on the indices) consists of $A=I_K$. Then by substituting this matrix it is easy to verify that the minimum point for the other parameters $\Lambda_k$ and $\pi_k$ occur at their true values $\tilde{\Lambda}_k$ and $\tilde{\pi}_k$, for $k = 1 ,\dots , K$. Note that for $\pi_i 's$ a lagrange multilplier method must be applied.
Since according to Theorem \ref{newmfvba}, $m=\mathfrak{M} (B)$ is the solution obtained from the equation MFVBA, 
therefore, a necessary and sufficient condition for the solution to the system \eqref{mfvba1} and \eqref{mfvba2} to be the true answer is that the solution $B$ of the system \eqref{bonyadi} coincide with $I_k$ upto permutation of the rows.

\begin{theorem}\label{theotheo}
The solution to the system of equations \eqref{bonyadi} gives rise to the true parameters $\{\tilde{\mu}_k, \tilde{\Lambda}_k, \tilde{\pi}_k\}_{k=1,\dots,K}$
if and only if $B$ coincides with one of the vertices of the simplex $(\Delta)^K$ where $\Delta=\{(x_1,\dots,x_{K})|\sum_i x_i =1\}$. 
\end{theorem}

According to \eqref{phihatag} is concave and it can be seen by \eqref{sydef} that $\psi$ is convex. Thus by adjusting the value of $\beta$ 
in the definition of $\lambda$ in \eqref{epla} one can turn $A\rightarrow (\hat{\phi}_A (\mathfrak{M} (A) +\psi (A)) $ into a concave function so that all it minimum points occur at the vertices of $\Delta$. 
From concavity one can deduce that the only critical points of $A\rightarrow (\hat{\phi}_A (\mathfrak{M} (A) +\psi (A))$ correspond to its absolute maximum. We claim that this critical point differs from $B$ given by Theorems \ref{the9} and \ref{newmfvba}.
In fact, at an absolute maximum point $B_0$ the Hessian of $A\rightarrow (\hat{\phi}_A (\mathfrak{M} (A) +\psi (A))$ will be non-positive definite while from \eqref{124} and \eqref{126}, it can be proved that The Hessian at $B$ must be positive definite. This can be proved by taking a path $t\rightarrow A_t$ with $A_0=B$ and compute $\frac{d^2}{dt^2} \left ( (\hat{\phi}_{A_t} (\mathfrak{M} (A_t) +\psi (A_t)) \right )$ by using and taking into account that $\mathfrak{M}(B)$ is a critical point of $\xi \rightarrow \hat{\phi}_B (\xi)$.

\begin{theorem}{PT}
By adjusting the temperature parameter $\beta$ one can assure that the solution to MFVI for $N$ large enough approaches the correct value of the parameters for $N$ large enough.
\end{theorem}

\section{Conclusion} 
In many high dimensional data science mining, one of the main modeling is the Gaussian Mixture Model (GMM). Direct exact computation with this model is very costly in computations. The Mean Field Variational Bayesian Inference (MFVBI) is classically used for approximate, but fast computation of the posterior probability density function (pdf) within this model. However, even if many properties of this model and this approximation computation are well-known, it suffers from lack of uncertainty quantification. 
In this paper, we forge foundation for a mathematical treatment of the MFVBI applied to the GMM. Several fundamental concepts from statistical mechanics such as partition function, Legendre transform, free energy and phase transition are revisited throughout our analysis. We can consider the GMM model as a generalization of field theory setup described very briefly and in an extremely non-rigorous manner in \cite{opper} for the Ising model. 
We have developed in this paper theoretical basis which elucidates towards which points the solutions to MFVI can converge and when the resulting point correspond to the correct values of the parameters. Furthermore, the temperature parameter included accommodates a simple modification whose adjustment can guarantee the accuracy 
of the mean field variational method. We defer the computational treatment of the problem to another paper. 
Finally it might be of worth to mention that the setup emerging throughout this paper seems to put forward new perspectives on a geometric understanding of statistical mechanics itself. 

%\vspace{1cm}

%\section*{Statements and Declarations}
%The authors did not receive support from any organization for the submitted work. 
%\\%
%On behalf of all authors, the corresponding author states that there is no conflict of interest. 
%\\
%The authors have no relevant financial or non-financial interests to disclose.

\appendix
\section*{Appendix} 
\section{Average computation}\label{app4} 
For $z\in \mathscr{Z}_{A_d}$ based on the definition of the true probability distribution $\tilde{p}_{N}$ in \eqref{ptrue}
\begin{align}\label{znkcomp2}
&\left\langle \sum_{n,k} z_{nk} (\mu_k- x_n)^T\Lambda_{k}( \mu_k- x_n)\right\rangle _{\tilde{p}_{_{N}}}\notag\\
&\quad\quad\quad=\left\langle\sum_{n, k',k}\tilde{z}_{nk'}z_{nk} \left((\mu_k-\tilde{\mu}_{k'}) - (x_n-\tilde{\mu}_{k'})\right)^T\Lambda_{k}\left((\mu_k-\mu_{k'})- (x_n-\tilde{\mu}_{k'})\right) \right\rangle _{\tilde{p}_{_{N}} (|)}\notag\\
&\quad\quad\quad =\sum_{k,k'} A_{z}(k,k')\left\langle ((\mu_k-\tilde{\mu}_{k'}) - (y-\tilde{\mu}_{k'}))^T\Lambda_{k}((\mu_k-\tilde{\mu}_{k'})- (y-\tilde{\mu}_{k'})) \right\rangle _{\mathcal{N}(y|\tilde{\mu}_{k'}, \tilde{\Lambda}_{k'})}\notag\\
& \quad\quad\quad=\sum_{k,k'} A_{z}(k,k')\sum_{i,j} \tilde\Lambda^{ij} _{k'} \Lambda_{k,ij}+\sum_{k,k'} A_{z}(k,k') (\mu_k -\tilde{\mu}_{k'})^T\Lambda_{k} (\mu_k -\tilde{\mu}_{k'})\notag\\
& \quad\quad\quad= \sum_{k,k'} \tilde{N}_{k'} \alpha_N (k,k')\sum_{i,j} \tilde\Lambda^{ij} _{k'} \Lambda_{k,ij}+\sum_{k,k'}\tilde{N}_{k'} \alpha_N (k,k') (\mu_k -\tilde{\mu}_{k'})^T\Lambda_{k} (\mu_k -\tilde{\mu}_{k'})
\end{align}

Here we are using relation \eqref{alfa1}. Likewise, $\tilde{p}_N(|)$ is the true probability measure \eqref{ptrue} conditioning prior knowledge about the true classes of the data. 
Similarly, we get
\begin{align}\label{znkcomp3}
\left\langle \sum_{k,n}\left( z_{nk} \log \pi _k 
+ \frac{1}{2}z_{nk} \log |\Lambda _k | \right) \right\rangle _{\tilde{p}_{N} (|)} \!\!\!\!= \sum_k \sum_{k'} \alpha_{N}(k,k')\tilde{N}_{k'}\left( \log \pi _k 
+ \frac{1}{2} \log |\Lambda _k | \right).
\end{align}
Let $h$ be a symmetric real valued function. A $U$-statistics with kernel $h$ of degree $m$ is 
\[
U_n = { \binom{n}{m}}^{-1} \sum\limits_{(i_1,\dots,i_m)\in I_{n,m}} h(X_{i_1},\dots,X_{i_m}),
\]
where 
$I_{n,m}=\left\{ (i_1,\dots,i_m)\ \ | 1\leq i_1<\dots<i_m\leq n\right\}$ is all possible ordered $m$-tuples. 
\begin{theorem}\cite[Theorem 1]{Chen}\label{chen}
If 
$E|h(X_{i_1},\dots,X_{i_m})|<\infty,$
then $U_n\rightarrow \theta(F)$ almost surly, 
where $F$ is common distribution function and $\theta$ represent the population mean (it's the parameter that the $U$-statistic aims to estimate).
\end{theorem}
\begin{theorem}\label{prop6}

The function 
\begin{equation}
\Phi (\Lambda, \mu ,\pi, z)=\sum_{k=1} ^K \bigg ( -\frac{1}{2\lambda_0}\sum_{n=1}^N z_{nk} (\nu_k- x_n)^T\Lambda_k (\mu_k - x_n)+\frac{1}{\lambda_0} \sum_{n}\left( z_{nk} \log \pi _k 
+ \frac{1}{2}z_{nk} \log |\Lambda _k | \right) \bigg )
\end{equation}
restricted to $M\times \mathscr{Z}_{A_d}$ converges with $N$, 
almost surely with respect to the distribution of the data $\{x_n\}_{ n=1}^N$, and uniformly over compact subsets of $M$ towards its average which according to the relations \eqref{znkcomp2} and \eqref{znkcomp3} is given by:
\begin{equation}\label{theo154}
\begin{split}
&\left\langle \Phi |_{M\times \mathscr{Z}_{A_d}} \right\rangle_{\tilde{p}_{N} }=\\
&\ - \sum_{k'} \frac{\tilde{N}_{k'}}{\lambda_0} \left(\frac{1}{2}\sum_{k} \alpha_N (k,k') \sum_{i,j} \tilde\Lambda^{ij} _{k'} \Lambda_{k, ij}+\frac{1}{2}\sum_{k} \alpha_N (k,k') \tilde{\pi}_{k'} (\mu_k -\mu_{k'})^T\Lambda_{k} (\mu_k -\mu_{k'})\right)\\
& \ + \sum_{k}\alpha_N (k,k')\left( \log \pi _k 
+ \frac{1}{2} \log |\Lambda _k | \right) \bigg ) 
\end{split}
\end{equation}
Here $\Lambda= (\Lambda_1,\dots,\Lambda_K), \mu =(\mu_1,\dots,\mu_K)$ and $\pi = (\pi_1,\dots,\pi_K)$.
\end{theorem}

\end{document}